[1994/12/01]
\documentclass[draft]{article}
\title{Solutions of Quasi-Geostrophic turbulence in multi-layered configurations}
\author{S. Jamal \\
School of Mathematics and Centre for Differential Equations,\\ Continuum Mechanics and Applications,\\
University of the Witwatersrand, Johannesburg, \\
Wits 2001, South Africa\\ Sameerah.Jamal@wits.ac.za}
\date{}

\usepackage{amsmath,amsthm}
\usepackage[final]{graphicx}
\usepackage{subcaption}

\chardef\bslash=`\\ 

\theoremstyle{definition}

\theoremstyle{remark}

\def\({\left(}
\def\){\right  )}
\def\dsy{\displaystyle}
\def\p{\partial}
\def\[{\left[}
\def\]{\right]}
\def\bn{\begin{equation}}
\def\en{\end{equation}}
\def\de{\delta}
\def\a{\alpha}

\newcommand{\eval}[2][\right]{\relax
  \ifx#1\right\relax \left.\fi#2#1\rvert}

\begin{document}
\maketitle
\markboth{Quasi-Geostrophic turbulence models}
{Quasi-Geostrophic turbulence models}
\renewcommand{\sectionmark}[1]{}

\begin{abstract}
We consider quasi-geostrophic (Q-G) models in two- and three-layers that are useful in theoretical studies of  planetary atmospheres and oceans.  In these models, the streamfunctions are given by  (1+2)  partial differential systems of evolution equations. A two-layer Q-G model, in a simplified version, is dependent exclusively on the Rossby radius of deformation.
However,  the $f$-plane  Q-G point vortex model  contains factors such as  the density,  thickness of each layer, the Coriolis parameter, and the constant of gravitational acceleration,  and this two-layered model admits a lesser number of  Lie point symmetries, as compared to the simplified model. Finally, we study a three-layer oceanography Q-G model of special interest,   which includes   asymmetric wind curl forcing or Ekman pumping, that drives double-gyre ocean circulation. In three-layers, we obtain solutions pertaining to the wind-driven double-gyre ocean flow for a range of physically relevant features, such as lateral friction and the analogue parameters of the $f$-plane  Q-G model. Zero-order invariants are used to reduce the  partial differential systems to ordinary differential systems. We determine conservation laws for these Q-G systems via multiplier methods. 
\end{abstract}

                             

\section{\label{sec:level1}Introduction}
Q-G theoretical models are of significance in meteorological and oceanographic analysis, notably because they describe  salient physical features, and have shown reasonable progress under numerical investigations. Q-G theory itself, contains robust observations that on large space and time scales, the evolution of the horizontal velocity and pressure gradient fields almost preserve a ``geostrophic'' balance between Coriolis and pressure gradients forces  \cite{con}. 
 
 In recent years, numerous investigations have revolved around  Q-G equations, in many different contexts.  For example, concerning the barotropic and two-layer Q-G models, some important results concerning the bifurcation behaviour and internal modes of variability have been obtained; see \cite{nauw} for an interesting discussion. The work by  \cite{jam} formulated a two-layer Q-G potential vortex  model as an infinite-dimensional Hamiltonian system and applied NoetherÕs theorem to derive conservation laws. Theoretical and experimental work has also been devoted to the study of single-layered atmospheric pressure models in the hopes of gaining insight into multi-layered   continuously stratified meteorological equations. Processes such as  baroclinic instability are well defined through research about two-layer atmospheric models.

In this paper, we begin by reviewing some of the known features of two- and three-layer non-dimensionalized  models in this area. Detailed discussions encompassing different points of theoretical views, as well as derivations of these models can be found in  \cite{ped,pan,sal}. 
In the literature pertaining to versions of Q-G  models,  the  velocity profile $\textbf{v}$ is obtained by standard techniques     from the Biot-Savart
law \cite{newton}
  $$\mathbf{v}=\left(\dot{x}, \dot{y}\right)=\bigtriangledown^{\perp}\psi=\left(-\frac{\p \psi}{\p x_2},\frac{\p \psi}{\p x_1}\right),$$
combined with an evolution equation for the geostrophic potential vorticity $Q$,
$$\frac{\p Q}{\p t}+\[\psi,Q\]=F(\psi)+D(\psi).$$
The Jacobian is given by $\dsy{[a,b]=\frac{\p a}{\p x_1}\frac{\p b}{\p x_2}-\frac{\p a}{\p x_2}\frac{\p b}{\p x_1},}$
$\(x_1,x_2\)$ are horizontal coordinates, $F$ and $D$ are forcing and
dissipation terms respectively, and $Q$ is defined by a linear differential operator $L$ in space variables
$Q=L(\psi).$
The functional form of  $L$ may be altered to provide different versions of the model. A generic equation is expressed by
$$L(\psi)=\triangle \psi+\beta_0x_2+p(x_3)\frac{\p}{\p x_3}\(q(x_3)\frac{\p \psi}{\p x_3}\),$$
where $\beta_0$  ($\ge0$) is a constant, $\triangle=\(\frac{\p}{\p x^2_{1}},\frac{\p}{\p x^2_{2}}\), $ $p(x_3)$ and  $q(x_3)$
are functions related to a reference state density structure \cite{con}.  
  The quantity $\dsy{\frac{\p \psi}{\p x_3}}$ represents temperature or buoyancy.  On the rapidly growing literature, there exists a hierarchy of models. For example, theoretical models that assume $\dsy{\frac{\p \psi}{\p x_3}=0}$, are called barotropic, and models that assume the contrary are called baroclinic. It is well known that barotropic versions with $\beta_0 = 0$ are  2D incompressible Navier-Stokes equations.

In analysing the physics of the atmosphere and ocean, layered models are of greater importance. 
The aforementioned geophysical fluid dynamics are modelled by sophisticated nonlinear partial differential equations (PDEs) whose solutions are derived predominantly from  numerical techniques. 
 Q-G models with  various layers are analysed,  by using the method of group invariant transformations, specifically the Lie point symmetries of the systems. The fundamental importance of Lie symmetries is that they reduce the number of independent variables in a DE and hence facilitate the solution. Moreover, group invariant transformations enable the classification of equations and provide information on DEs, and subsequently for the models.
Therefore, based on the discussion presented above, the objective of this study is to document the Lie  point symmetries admitted, in order to obtain solutions for the two- and three-layer Q-G problem.

 Also, in this paper, we  utilise the multiplier approach which was introduced by Stuedel \cite{ste}. This method relies on the well known result that the Euler-Lagrange operator annihilates a total divergence \cite{here}. The role of multipliers  has been shown to play a significant role in the construction of conserved quantities \cite{euler} and in determining hierarchies of PDEs. Essentially, a knowledge of a multiplier,  leads to a conserved flow. To calculate the conserved densities and fluxes, one requires the inversion of the total divergence operator, and hence the usage of a homotopy operator.  The homotopy operator \cite{o,ab,here} is a powerful algorithmic tool that originates from homological algebra and variational bi-complexes. 
 
 To our knowledge, there are no comprehensive formulations of the Lie symmetries and conservation laws derived from multipliers,  of these geophysical systems and we provide this,  as well as  present some solutions associated with configurations of two-layers, and a three-layer oceanography model. 
The paper is organized as follows. The next section reviews  the mathematical analysis required
to perform our investigation. First, in section 3, we discuss three important Q-G models and their  Lie point symmetries to determine invariant solutions.  For the Rossby radius deformation model and ocean three-layer model, the closed-form solutions are expressed in terms of sine and cosine functions, whereas for the $f-$plane two-layer model, the closed-form solution is expressed in terms of exponential functions. 
Section 4 describes the multipliers of two of the Q-G models under study and the corresponding conservation laws.
 Conclusions are given in section 5.


\setcounter{equation}{0}

\section{Mathematical Description}
\setcounter{equation}{0}
In this section, we outline the general procedure for determining point symmetries for an arbitrary system of equations.  Consider a nonlinear system with  $q$  unknown functions $u^a$ which depends on $p$ independent variables  $x^i$, i.e. we denote
$u = (u^1,\ldots, u^q)$ and $x= (x^1,\ldots, x^p)$, respectively. 
Let \bn\label{main} G_\alpha\(x,u^{(k)}\)=0, \quad \alpha=1, \ldots, q,\en
be a system of $m$ nonlinear differential equations, where $u^{(k)}$ represents the $k^{th}$  derivative of $u$ with respect to $x$. 
A one-parameter Lie group of transformations ($\epsilon$ is the group parameter) that is invariant under (\ref{main}) is given by \bn\label{main2} \bar{x}=\Xi (x,u;\epsilon)\quad \bar{u}=\Phi (x,u;\epsilon).\en
Invariance of (\ref{main})  under the transformation (\ref{main2})  implies that any solution $u = \Theta(x)$ of  (\ref{main}) maps into another solution $v = \Psi(x; \epsilon) $ of  (\ref{main}).  
 Expanding (\ref{main2}) around the identity $\epsilon = 0$, we can generate the following infinitesimal transformations:
\bn\begin{array}{ll}
&\bar{x}^i=x^i+\epsilon\xi^i(x,u)+{\cal O}(\epsilon^2),\quad i=1, \ldots, p,\\
&\bar{u}^\alpha=u^\alpha+\epsilon\eta^\alpha(x,u)+{\cal O}(\epsilon^2). 
\end{array}
\en
The action of the Lie group can be recovered from that of its
infinitesimal generators acting on the space of independent and dependent variables. Hence, we consider the following infinitesimal vector field
\bn\label{xop}\dsy{X=\xi^i\p_{x^i}+\eta^\alpha\p_{u^{\alpha.}}}\en
The action of X is extended to all derivatives appearing in the equation in question through the appropriate prolongation. The infinitesimal criterion for invariance is given by \bn\label{crib}X\[\textrm{LHS Eq.} (\ref{main}) \]\mid_{Eq. (\ref{main}) }=0.\en Eq.  (\ref{crib}) yields an overdetermined system of linear homogeneous equation which can be solved algorithmically, more details can be found in  \cite{o} among other texts.
The total differentiation operator $D_i$ with respect to $x^i$ is given by
\bn \label{eq1}
D_i={{\p}\over{\p x^i}}+u_i^\alpha {{\p}\over{\p u^\alpha}}+
u_{ij}^\alpha {{\p}\over{\p u_j^\alpha }} + \ldots~~~
\en
The Euler operator is given by
\bn {{\de}\over{\de
u^\alpha}}={{\p}\over{\p u^\alpha}}+\sum_{s\geq1}(-1)^sD_{i_1}\cdots
D_{i_s}\;{{\p}\over{\p u^\alpha_{i_1\cdots i_s}}}.
\label{3.12} \en
A current ${\bf T}=(T^1,\ldots,T^n)$ is conserved if it satisfies
\bn\label{con-law}
D_i \,  T^i=0
\en
along the solutions of (\ref{main}), and Eq. (\ref{con-law}) is called a local conservation law.
The operator in Eq. (\ref{xop}) can be used to define the Lagrange system
$$\frac{dx^i}{\xi^i}=\frac{du^\alpha}{\eta^\alpha}$$
whose solution provides the zero-order invariants
\bn\label{in}W^{[0]}(x^i,u^\alpha).\en
These invariants can be used in order to reduce the order of the PDE. Further details of the relevant equations and formulae can be found in, inter alia, \cite{hs}. 
The multipliers $\Lambda^\alpha$ of (\ref{main}) satisfies the relation
\bn\label{con}
D_i \,  T^i=\Lambda^\alpha G_\alpha
\en
for the function $u^\alpha$ \cite{o,ste}. The determining equations for the multipliers are then
\bn\label{e}{{\de}\over{\de
u^\alpha}} [\Lambda^\alpha G_\alpha]=0.\en
Eq. (\ref{con}) is satisfied for the arbitrary functions $u^\alpha$ and not only for the solutions of (\ref{main}). Conserved vectors may be derived systematically using (\ref{con}) as the determining equation, however, in some cases it is simple to construct the conserved vectors by elementary manipulations once the multiplier has been determined. 
Alternatively, the conserved quantities are determined by a homotopy operator, see \cite{here} for details. 



\section{Layered Q-G models}
Below, we perform a symmetry analysis and  determine the group invariant solutions. In particular, we focus our analysis on two models of special interest which consist of two layers: the $f$-plane  Q-G point vortex model and the simplified model dependent only on the Rossby radius of deformation. Moreover, we study a three-layer oceanography Q-G model. The latter is a model which includes   asymmetric wind curl forcing (Ekman pumping) that drives double-gyre ocean circulation.

\subsection{Model I: Two-layer Q-G point vortex model}

The basic evolution equations for the $f$-plane two-layer Q-G potential vortex model are represented by the evolution equations  \cite{ped,jam}
\bn\label{q1} 
\dsy{\frac{\p \omega_i }{\p t} +[\omega_i,\psi_i]=0, \quad i=1,2,}
\en
where the potential vorticity $\omega_i$  is expressed as \bn\label{w1}\dsy{\omega_i=\triangle \psi_i+\epsilon_i (\psi_{2}-\alpha_i\psi_1) }, \quad \textrm{and}\en
$$\a_1=1, \quad \a_2=\frac{\rho_1}{\rho_2},\quad \epsilon_i=\frac{l^2\rho_i}{(\rho_{2}-\rho_1)gH_i}.$$
In this model,  $\psi(x,y,t)$ is the corresponding streamfunction in the $i$th layer, $\rho_i$ is the respective fluid densities $\rho_i<\rho_{i+1}$, $H_i$ is the thicknesses of each layer, $l$ the Coriolis parameter, and $g$ the constant of gravitational acceleration. The thickness of each layer is assumed to be much smaller than the horizontal scale of the system so that the fluid is always in hydrostatic equilibrium,  the upper boundary is free and the lower boundary is solid \cite{jam}. Before we proceed with the symmetry analysis, we remark that  the system of equations (\ref{w1}) are Poisson equations with the vorticity on the left-hand side. Exercising the symmetry condition Eq. (\ref{crib}), we obtain a system of twenty-five equations  (for the derivation of this system, we used Maple software). Hence, Eq. (\ref{q1})   admits a Lie algebra of point symmetries that   is 7-dimensional, viz.,  
\bn \begin{array}{ll}
&X_1=\p_x,\quad X_2=\p_y,\quad X_3=\p_t,\quad X_4=\p_{\psi_1},\quad X_5=\p_{\psi_2},\\
&X_6=y\p_x-x\p_y,\quad X_7=t\p_t-\psi_1\p_{\psi_1}-\psi_2\p_{\psi_2}.
\label{ssa}\end{array}\en
The commutator of symmetry generators are defined, in general, as $\[X_i,X_j\]=X_iX_j-X_jX_i$. Non-vanishing commutator relations in this case are $\[X_1,X_6\]=-X_2,\, \[X_2,X_6\]=X_1,\, \[X_3,X_7\]=X_3,\, \[X_4,X_7\]=-X_4,\, \[X_5,X_7\]=-X_5. $
In the following subsection, we consider a special two-layered model  which has been proposed in the literature for its simplicity.

\subsection{Model II: A simplified two-layered model}
A second two-layer model that is simpler but shares many of the features of the full two-layer model above, has potential vorticity  expressed as
\bn\label{w2}\omega_i=\triangle \psi_i+(-1)^i\lambda^{-2} (\psi_{1}-\psi_{2}) ,\en with Eq. (\ref{q1}) unchanged.
This system contains only one parameter, the Rossby radius of deformation $\lambda$, which can be interpreted as the rigidity of the interface separating the two layers \cite{jam}. Infinitesimal generators of every one parameter Lie group of point symmetries, from the Lie symmetry condition Eq. (\ref{crib}),  are ($a(t),\, b(t),\, c(t)$ are arbitrary):
\bn \begin{array}{ll}
&Y_1=\p_t,\quad Y_2=\p_{\psi_2},\quad Y_3=y\p_x-x\p_y,\quad Y_4=t\p_t-\psi_1\p_{\psi_1}-\psi_2\p_{\psi_2},\\
&Y_5=a(t)\[\p_{\psi_1}+\p_{\psi_2}\],\quad Y_6=b(t)\p_x-yb'(t)\p_{\psi_1}-yb'(t)\p_{\psi_2},\\
& Y_7=c(t)\p_y+xc'(t)\p_{\psi_1}+xc'(t)\p_{\psi_2},\\
&Y_8={ty\p_x-tx\p_y-\frac12 (x^2+y^2)\p_{\psi_1}-\frac12 (x^2+y^2)\p_{\psi_2}}.
\label{ssa}\end{array}\en
A variant of this model has been analysed in terms of its Lie symmetries \cite{b1}. Next, we look at the Lie symmetries of a three-layered model.

\subsection{Model III: Three-layer Q-G ocean circulation model}
Oceans embody  the largest and least understood elements of our extensive climate system. In this analysis, we adopt the model dynamics of the Q-G system described by  \cite{nauw}, where the authors focussed on  flow regimes  connected to ocean circulations. We consider a double-gyre three-layer Q-G ocean model  describing idealised midlatitude ocean circulation. Double-gyre flow  in a rectangular basin contributes to the understanding of the variability in the upper ocean \cite{jia} and
Q-G models, in essence,  aid in  understanding the mechanism of the double-gyre phenomenon. The system is given by
\bn\begin{array}{ll}\label{q3} 
&\dsy{\frac{\p \omega_1 }{\p t} +[\omega_1,\psi_1]=A_H \bigtriangledown^4 \psi_1 +\frac{f_0}{H_1} \Omega_e}\\
&\dsy{\frac{\p \omega_i }{\p t} +[\omega_i,\psi_i]=A_H \bigtriangledown^4 \psi_i\quad i=2,3,}
\end{array}\en
where the potential vorticities in each layer are  expressed as  \bn\begin{array}{ll}\label{w3}\dsy{\omega_1}&\dsy{=\bigtriangledown^2\psi_1-\frac{f_0}{H_1}h_1 +\beta_0y},\\
\dsy{\omega_2}&\dsy{=\bigtriangledown^2 \psi_2-\frac{f_0}{H_1} (h_2-h_1)+\beta_0y,}\\
\dsy{\omega_3}&\dsy{=\bigtriangledown^2 \psi_3+\frac{f_0}{H_3}h_2 +\beta_0y.}
\end{array}\en
The stratification within this model is characterized by three stacked layers of water with constant densities $\rho_i<\rho_{i+1}$ ($ i=1,2,3$) and mean layer thicknesses $H_i$ where $H=\sum_i H_i$. Lateral friction with a constant lateral friction coefficient $A_H$ is considered as the only dissipation mechanism and the domain is a square basin with dimension $2L\times 2L$. The Coriolis parameter, $f_0$, is taken at  $\theta_0$, which is the location of the mid-axis of the basin. The Coriolis parameter varies linearly with latitude as $f(y)=  f_0+\beta_0y$, where $\beta_0$ represents its meridional gradient. The quantities $h_i$ represent interface perturbations given by $$h_i= \frac{  f_0\(\psi_i    -\psi_{i+1}\)}{g'_i},\quad i=1,2.$$  The reduced gravity parameters are  $$g'_i= \frac{\(\rho_{i+ 1}-\rho_i\)}{\rho_0}\quad i=1,2.$$
The surface wind stress,  $\bf{\tau}=(\tau^x,\tau^y)$ is given by the analytic function
$$\tau^x=\tau_0\[\(1-2\alpha\frac{y}{L}\)\cos\(\frac{\pi y}{L}\)+\frac{2\alpha}{\pi}\sin\(\frac{\pi y}{L}\)\],\quad\tau^y=0$$
where $\tau_0$ is the amplitude. The asymmetric wind curl forcing (Ekman pumping) drives the double-gyre ocean circulation, and is given by
$$\Omega_e=\mu_0\(1-2\alpha\frac{y}{L}\)\sin\(\frac{\pi y}{L}\),\quad y\in[-L,L]$$
 and $\dsy{\mu_0=\frac{\tau_0\pi}{\rho_0f_0L}}$. The parameter $\alpha$  is introduced to control the amount of asymmetry in the Ekman pumping with respect to the mid-axis of the basin $( y=0)$. The Lie group of point symmetries of this Q-G system is 5-dimensional, generated by the vector fields  ($f(t)$ is arbitrary):
\bn \begin{array}{ll}
&Z_1=\p_t,\quad Z_2=\p_x,\quad Z_3=\p_{\psi_2},\quad Z_4=\p_{\psi_3},\\
&Z_5=f(t)\[\p_{\psi_1}+\p_{\psi_2}+\p_{\psi_3}\].
\label{ssa}\end{array}\en
The  numerical implementation used here is originally described in \cite{land} and, thereafter in \cite{nauw}.
\begin{table}[h!t]%
\caption{Model parameters for the three-layer Q-G system}
\centering
\begin{tabular}{|c|c|c|c|}
\hline
 \textbf{Parameter} &  \textbf{Value} & \textbf{Parameter} &  \textbf{Value}\cr
\hline\hline
$f_0$ & $1.0\times10^{-4}\textrm{s}^{-1}$ &$H$ & $4.0\times10^{3}\textrm{m}$  \cr
$\beta_0$ & $1.6\times10^{-11}\textrm{ms}^{-1}$ &$A_H$ & $3.0\times10^{2}\textrm{m}^{2}\textrm{s}^{-1}$ \cr
$\rho_0$ & $1.0\times10^{3}\textrm{kg }\textrm{m}^{-3}$ &$\mu_0$ & $2.5\times10^{-6}\textrm{m}\textrm{s}^{-1}$ \cr
$H_1$ & $6\times10^{2}\textrm{m}$ &$g'_1$ & $2.0\times10^{-2}\textrm{m}\textrm{s}^{-2}$  \cr
$H_2$ & $1.4\times10^{3}\textrm{m}$ &$g'_2$ & $3.0\times10^{-2}\textrm{m}\textrm{s}^{-2}$  \cr
$H_3$ & $2.0\times10^{3}\textrm{m}$ &$\alpha$ &$0.0$\cr
$\tau_0$ & $1.0\times 10^{-1}\textrm{N}\textrm{m}^{-2}$ & &\cr
\hline\hline
\end{tabular}
\label{table7}
\end{table}

\subsection{Analytic Solutions}

\begin{figure}
\begin{subfigure}{.4\textwidth}
  \centering
  \includegraphics[width=.8\linewidth]{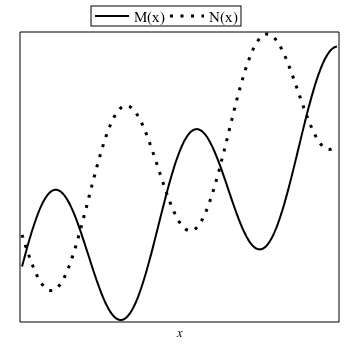}
  \caption{$\lambda=10^{-1}$}
  \label{fig:sfig1}
\end{subfigure}%
\begin{subfigure}{.4\textwidth}
  \centering
  \includegraphics[width=.8\linewidth]{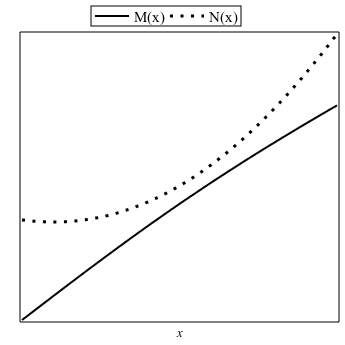}
  \caption{$\lambda=1$}
  \label{fig:sfig2}
\end{subfigure}
\begin{subfigure}{.4\textwidth}
  \centering
  \includegraphics[width=.8\linewidth]{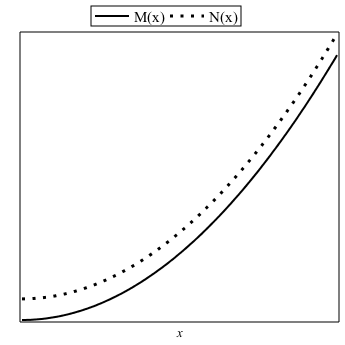}
  \caption{$\lambda=10^{2}$}
  \label{fig:sfig3}
\end{subfigure}%
\begin{subfigure}{.4\textwidth}
  \centering
  \includegraphics[width=.8\linewidth]{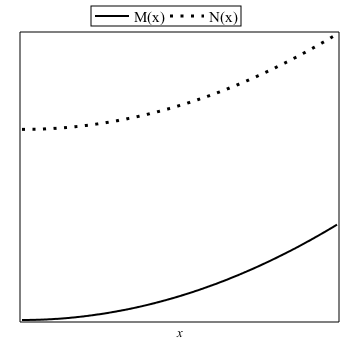}
  \caption{$\lambda=10^{3}$}
  \label{fig:sfig4}
\end{subfigure}%
\newline
\begin{subfigure}{.4\textwidth}
  \centering
  \includegraphics[width=.8\linewidth]{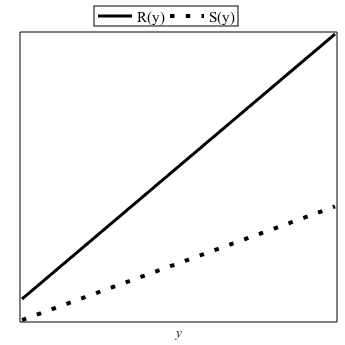}
  \caption{Solutions to Eqs. (\ref{pic1})}
  \label{fig:sfig5}
\end{subfigure}%
\begin{subfigure}{.4\textwidth}
  \centering
  \includegraphics[width=.8\linewidth]{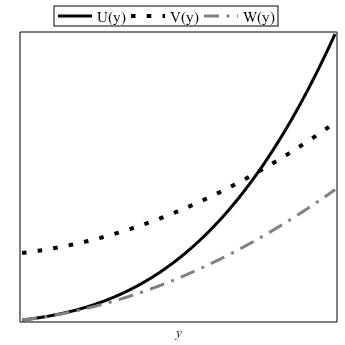}
  \caption{Solutions to Eqs. (\ref{ode4})}
  \label{fig:sfig6}
\end{subfigure}%
\caption{Progression of the solutions  are depicted. Typical values for the Rossby radius $\lambda$  in Eqs. (\ref{pic2})  are selected. In these numerical solutions, we choose for (a)-(d) $\sigma_1=\sigma_2=1;$ (e) $\mu_1=\mu_3=\mu_4=1, \mu_2=0.5$; (f) $\kappa_1=0.5, \kappa_2=1, \kappa_3=1.5.$ }
\label{fig:fig}
\end{figure}

The knowledge of Lie point symmetries of a DE makes possible the determination of solutions to these equations, which are invariant under a given point symmetry. 
In this subsection we apply the Lie symmetries and resultant  zero-order invariants,  to reduce the PDE systems under study and determine analytic solutions. We discuss a reduction for each model, and display the evolution of some solutions of interest in Figure  \ref{fig:fig}. Moreover, our models are (1+2) evolution equations and, in order to reduce  to an ordinary differential equation (ODE) system, we have to apply the invariants of two Lie point symmetries. Often, after the application of one symmetry, one obtains a reduced equation which in turn admits reduced symmetries. Reduced symmetries are determined in the standard way. 

\textbf{Model I:}
For the $f$-plane two-layer Q-G system, as we saw above, Eq. (\ref{q1}) admits a 6-dimensional Lie algebra of point symmetries.  Therefore, suppose we perform reductions, successively with the  linear combinations  ${\cal H}_\alpha=X_3+\mu_1X_4+\mu_2X_5$ followed by ${\cal H}_\beta=X_1+\mu_3X_4+\mu_4X_5.$
Hence, applying the group invariants of the sub-algebra $\Sigma_I=\{{\cal H}_\alpha,{\cal H}_\beta\}$,  
 leads to invariant solutions, which are the streamfunctions,
\bn\begin{array}{ll}\psi_1(t,x,y)&=\mu_1t+\mu_3x+R(y),\\
 \psi_2(t,x,y)&=\mu_2t+\mu_4x+S(y),  \label{sol1}\end{array}\en
where $R(y)$ and $S(y)$ satisfy the ODE system
\bn\begin{array}{lcl}\label{pic1}
&\dsy{(\mu_1-\mu_2)l^2\rho_1-\mu_4l^2\rho_1R'+\mu_3l^2\rho_1S'-(\rho_1-\rho_2)gH_1\mu_3R^{'''}=0},\\
&\dsy{\mu_1l^2\rho_1-\mu_2l^2\rho_2-\mu_4l^2\rho_1R'+\mu_3l^2\rho_1S'-(\rho_1-\rho_2)gH_2\mu_4S^{'''}=0,}\\
\end{array}\en
and dashes denote differentiation with respect to $y$. We continue with the determination of the group invariant solutions for Model II.

\textbf{Model II:} 
Next, we derive solutions for the simpler two-layer model (\ref{w2}), and recall that this system admits 8 Lie point symmetries.  
The application of the linear combination $Y_1+Y_5$ followed by the reduced symmetry  $\dsy{\p_y+\sigma_1\psi_1+\sigma_2\psi_2}$ produces streamfunctions of the form 
\bn\begin{array}{ll}\psi_1(t,x,y)&=\int a(t) dt+ \sigma_1y +M(x),\\
\psi_2(t,x,y)&=\int a(t) dt+ \sigma_2y +N(x).
   \label{sol3}\end{array}\en
$M(x), N(x)$  satisfy the reduced system below and dashes denote differentiation with respect to $x$: 
\bn\begin{array}{ll}\label{pic2}
&\sigma_1N'-\sigma_2(M'+\lambda^2N^{'''})=0,\\
&-\sigma_1N'+\sigma_2M'-\sigma_1\lambda^2M^{'''}=0.\\
\end{array}\en
In Figures (\ref{fig:sfig1})-(\ref{fig:sfig4}), we graph the numerical evolution of $M(x)$ and $N(x)$ for various values of the parameter $\lambda$.

\textbf{Model III:}
Following the steps of the previous model, we find the invariant solution of the ocean model. Information contained in Table \ref{table7}, provides the parameter values.
Here we apply the linear combinations $Z_1+Z_5$ and the reduced generator $Z_2+\kappa_1\psi_1+\kappa_2\psi_2+\kappa_3\psi_3$, to find streamfunctions with functional form
\bn\begin{array}{ll}\psi_1(t,x,y)&=\int f(t) dt+ \kappa_1x +U(y),\\
\psi_2(t,x,y)&=\int f(t) dt+ \kappa_2x +V(y),\\
\psi_3(t,x,y)&=\int f(t) dt+ \kappa_3x +W(y).
\label{sol4}\end{array}\en
 $U(y),\,V(y),\, W(y)$ are given by the ODE system
\bn\begin{array}{ll}\label{ode4}
&\dsy{-f_0^2 L^2 \kappa_1 \rho^2_0 V' + \kappa_2 f_0^2 L^2 \rho^2_0 U' + (\rho_1 - \rho_2) \Big(\beta_0 H_1 L^2  \kappa_1  \rho_0 - L \pi \tau_0 \sin\(\frac{\pi y}{L}\) } \\
 &\phantom{mmmmmm}    \dsy{+2 \pi y \alpha \tau_0 \sin\(\frac{\pi y}{L}\) + 
    H_1 L^2 \kappa_1  \rho_0 U^{'''}- A_H H_1 L^2  \rho_0 U^{''''}\Big)=0,}\\
&\dsy{f_0^2 (\kappa_2 - \kappa_1) \rho_0 \(\rho_2 - \rho_3\)V'+ \kappa_2 \Big(f_0^2 (\rho_1 - \rho_2) \rho_0 W' + (\rho_2 - \rho_3) \Big(-f_0^2 \rho_0 U'}\\
&\phantom{mmmmmmmmmmmmmmmmmm}  \dsy{- H_2 (\rho_1 - \rho_2) \(\beta_0 + V^{'''}\)\Big)\Big)=0},\\
&\dsy{-H_2 f_0^2 \kappa_3 \rho_0 V' + f_0^2 \( \kappa_2 H_3 - H_3  \kappa_3 + H_2  \kappa_3\)  \rho_0 W'}\\
&\phantom{mmmmmmmmmmmmmmmmm}  \dsy{+ H_3 H_2 \kappa_3 \(\rho_2 - \rho_3\) \(\beta_0 +W^{'''}\)=0, }\end{array}\en
where dashes denote differentiation with respect to $y$. In Figure \ref{fig:fig}, we give numerical solutions of the ODE systems.

\section{Conservation Laws of Model I and II}
A vital counterpart of any symmetry study  of a PDE is information about associated conservation laws. It is well known that these conservation laws play an important role in mathematical physics as they describe essential physical properties of the modelled process. Conservation laws may be used in the numerical integration of PDEs, for example, to control numerical errors \cite{rj}. The multiplier approach is a well known technique, used to derive conservation laws.  We study the two cases, Model I and II, which are the two-layered Q-G systems.

To implement this approach, consider a multiplier comprising of the independent variables, the dependent variable  and derivatives of dependent variables up to some fixed order. Multipliers $\Lambda^i\,(i=1,2)$, of the system of equations (\ref{q1})
have the property that $$\sum_{i=1}^2\Lambda^i\[\textrm{LHS of Eqs.}  (\ref{q1}) \]=D_xT^x+D_yT^y+D_tT^t$$
for all streamfunctions $\psi_i(x,y,t)$, where the total derivative operative is defined as in Eq. (\ref{eq1}). The right-hand side of this expression is a divergence expression and $T^j\; (j=x,y,t)$ are the components of the conserved vector
$T=(T^x,T^y,T^t)$. 
The determining equations for the multipliers $\Lambda^i$ are obtained from the expressions
\bn\label{el}{{\de}\over{\de
\psi_1}}\[\sum_{i=1}^2\Lambda^i  (\textrm{LHS Eq.} \ref{q1})\]=0,\quad {{\de}\over{\de
\psi_2}}\[\sum_{i=1}^2\Lambda^i  (\textrm{LHS Eq.} \ref{q1})\]=0, \en
where $ \dsy{{\de}\over{\de
\psi_i}}$  are the  Euler operators defined in (\ref{3.12}), which annihilate divergence expressions. This computation is extensive, especially since the order of the multiplier increases the complexity of this step. We separate  Eq. (\ref{el}), and after expansion, with respect to different combinations of derivatives of $\psi_i$, we obtain an overdetermined system. We omit these tedious and lengthy calculations. Ultimately, we find the following multipliers of Eq. (\ref{q1}) in each case considered below. 

\textbf{Model I:} The form $\Lambda^i(x,y,t,\psi_1,\psi_2) $, provides the following multipliers; that is $\Lambda^i=\(\Lambda^1_I,\Lambda^2_I\)$, i.e.:
 \bn \begin{array}{cc}\label{a}
&\dsy{\(0,1\),\quad \(\psi_1,-\frac{H_2}{H_1}\psi_2\),\quad \(F_1(t),\frac{H_2}{H_1}F_1(t)\)},\\
&\dsy{\(F_2(t)x,-\frac{H_2}{H_1}x F_2(t)\),\quad \(F_3(t)y,-\frac{H_2}{H_1}yF_3(t)\)},\\
&\(\dsy{\frac12(x^2+y^2)F_4(t)},\dsy{-\frac{H_2}{H_1}\frac12(x^2+y^2)F_4(t)}\),\\
\end{array}\en
where $F_k\, (k=1,2,3,4)$ are arbitrary functions of $t$.

Each of the multipliers lead to a conserved vector. For instance, the components of the conserved vector  associated with the multiplier (0,1) above, are 
$$\begin{array}{lcl}
&T^t_I=&\dsy{\frac{l^2 \rho_1 \psi_1 -l^2 \rho_2 \psi_2 }{g H_2 \rho_1-g H_2 \rho_2},}\\
&T^x_I=&\dsy{
\frac{l^2 \rho_1 \psi_1 \frac{\p \psi_2}{\p y}}{g H_2 \rho_1-g H_2 \rho_2}-\psi_2  \frac{\p^3 \psi_2}{\p x^2\p y}+\frac{\p^2 \psi_2}{\p t \p x},} \\
&T^y_I=&\dsy{ -\frac{1}{g H_2 ( \rho_1 - \rho_2)} \Bigg(l^2 \rho_1\psi_1 \frac{\p \psi_2}{\p x}+g H_2 (\rho_1 - \rho_2)  \Big(\frac{\p^2 \psi_2}{\p y^2} \frac{\p \psi_2}{\p x}-\frac{\p \psi_2}{\p y}\frac{\p^2 \psi_2}{\p x\p y}}\\
& &\dsy{-\psi_2  
\frac{\p^3 \psi_2}{\p x^3}-\frac{\p^2 \psi_2}{\p t\p y} \Big) \Bigg).}
\end{array}$$

\textbf{Model II:} If we consider the multiplier $\Lambda^i\(x,y,t,\psi_1,\psi_2, \dsy{\frac{\p\psi_1}{\p x}}, \dsy{\frac{\p\psi_1}{\p y}}, \dsy{\frac{\p\psi_1}{\p t}}\)$ where we find six multipliers $\(\Lambda^1_{II},\Lambda^2_{II}\):$
 \bn \begin{array}{cc}\label{a2}
&\dsy{\(0,1\),\quad \(\psi_1,\psi_2\),\quad \(J_1(t),J_1(t)\),}\\
&\dsy{\(J_2(t)x, J_2(t)x\),\quad \(J_3(t)y,J_3(t)y\),}\\
&\(\dsy{\frac12(x^2+y^2)J_4(t)},\dsy{\frac12(-4\lambda^2 +x^2+y^2)J_4(t)}\),\\
\end{array}\en
where $J_k\, (k=1,2,3,4)$ are arbitrary functions of $t$.
Similarly, as an example, $\(J_1(t),J_1(t)\)$ yields the conserved vector
$$\begin{array}{ll}
T^t_{II}=&{\frac{J_1}{3}\( \frac{\p^2 \psi_1}{\p y^2}+ \frac{\p^2 \psi_1}{\p x^2}+\frac{\p^2 \psi_2}{\p y^2}+ \frac{\p^2 \psi_2}{\p x^2}\),}\\
T^x_{II}=&{\frac{1}{6 \lambda ^2}\Bigg(-2 \lambda ^2 J_1' \frac{\p \psi_1}{\p x}+J_1 \Bigg(-3 \psi_2  \frac{\p \psi_1}{\p y}+\psi_1  \(3 \frac{\p \psi_2}{\p x}+2 \lambda ^2 \(\frac{\p^3 \psi_1}{\p y^3}+\frac{\p^3 \psi_1}{\p x^2\p y}\)\)}\\
&{+4 \lambda ^2 \(\frac{\p \psi_1}{\p x} \frac{\p^2 \psi_1}{\p x\p y}-\frac{\p \psi_1}{\p y} \frac{\p^2 \psi_1}{\p x^2}+\frac{\p^2 \psi_1}{\p x\p t}\)\Bigg)\Bigg)+\frac{1}{6 \lambda ^2}\Bigg(-2 \lambda ^2 J_1' \frac{\p \psi_2}{\p x}+J_1}\\
& {\(-3 \psi_1  \frac{\p \psi_2}{\p y}+\psi_2  \Bigg(3 \frac{\p \psi_1}{\p y}+2 \lambda ^2 \(\frac{\p^3 \psi_2}{\p y^3}+\frac{\p^3 \psi_2}{\p x^2 \p y}\)\)+4 \lambda ^2 \Bigg(\frac{\p \psi_2}{\p x} \frac{\p^2 \psi_2}{\p x\p y}}\\
&{-\frac{\p \psi_2}{\p y} \frac{\p^2 \psi_2}{\p x^2}+\frac{\p^2 \psi_2}{\p t\p x}\Bigg)\Bigg)\Bigg)},\\
T^y_{II}=&{\frac{1}{6 \lambda ^2}\Bigg(-2 \lambda ^2 J_1' \frac{\p \psi_2}{\p y}+J_1 \Bigg(3 \psi_1 \frac{\p \psi_2}{\p x}-\psi_2  \(3 \frac{\p \psi_1}{\p x}+2 \lambda ^2 \(\frac{\p^3 \psi_2}{\p x\p y^2}+\frac{\p^3 \psi_2}{\p x^3}\)\)}\\
&{+4 \lambda ^2 \(\frac{\p^2 \psi_2}{\p y^2} \frac{\p \psi_2}{\p x}-\frac{\p \psi_2}{\p y} \frac{\p^2 \psi_2}{\p x\p y}+\frac{\p^2 \psi_2}{\p t \p y}\)\Bigg)\Bigg)+\frac{1}{6 \lambda ^2}\Bigg(-2 \lambda ^2 J_1' \frac{\p \psi_1}{\p y}+J_1}\\
&{ \Bigg(3 \psi_2  \frac{\p \psi_1}{\p x}+4 \lambda ^2 \frac{\p^2 \psi_1}{\p y^2}\frac{\p \psi_1}{\p x}-3 \psi_1 \frac{\p \psi_2}{\p x}-4 \lambda ^2 \frac{\p \psi_1}{\p y} \frac{\p^2 \psi_1}{\p x\p y}-2 \lambda ^2 \psi_1  \frac{\p^3 \psi_1}{\p x\p y^2}}\\
&{-2 \lambda ^2 \psi_1  \frac{\p^3 \psi_1}{\p x^3}+4 \lambda ^2 \frac{\p^2 \psi_1}{\p t\p y}\Bigg)\Bigg).}
\end{array}$$

\section{Conclusion}

In this study, we investigated algebraic properties or Lie symmetries, of the $f$-plane two-layer Q-G potential vortex model in hydrostatic equilibrium. Additionally, we considered a model dependent on the Rossby radius of deformation, and showed that this model is invariant under eight Lie point symmetries. Lastly, we investigated the symmetries and analytic solutions of a three-layer Q-G ocean circulation model. Our analysis led to the observation that the ocean model admits fewer symmetries than the two-layered models.

Moreover, the inherent zero-order invariants of the Lie symmetries were used to reduced the PDE models to ODE systems. The analogous symmetry reduction  of the simplified model, produced closed-form solutions with variable $\lambda$.
Concerning each model, we obtained an invariant solution and we illustrated  the progression of the model's solution. 
To this end, we employed multiplier methods to construct conservation laws for the two-layer Q-G systems.

\textbf{Acknowledgments}
The author would like to thank the National Research Foundation of South Africa for financial support, grant number 99279.


\end{document}